\documentclass[12pt,a4paper]{article}
\usepackage[english]{babel}
\usepackage{amsmath}
\usepackage{amssymb}
\usepackage{amsfonts}
\usepackage{hyperref}

\begin{document}


\renewcommand{\refname}{References}
\renewcommand{\contentsname}{Contents}

\begin{center}
{\huge Global unique solvability\\ of the initial-boundary value problem for the equations of\\ one-dimensional polytropic flows\\ of viscous compressible multifluids}
\end{center}

\medskip

\begin{center}
{\large Alexander Mamontov,\quad Dmitriy Prokudin}
\end{center}

\medskip

\begin{center}
{\large September 18, 2018}
\end{center}

\medskip

\begin{center}
{
Lavrentyev Institute of Hydrodynamics, \\ Siberian Branch of the Russian Academy of Sciences\\
pr. Lavrent'eva 15, Novosibirsk 630090, Russia}
\end{center}

\medskip

\begin{center}
{\bfseries Abstract}
\end{center}


\begin{center}
\begin{minipage}{110mm}
We consider the equations which describe polytropic one-dimensional flows of viscous compressible multifluids. We prove global existence and uniqueness of a solution to the initial-boundary value problem which corresponds to the flow in a bounded space domain.
\end{minipage}
\end{center}

\bigskip


\newpage

\tableofcontents

\bigskip

\section{Introduction}

\noindent\indent
The paper is devoted to the analysis of the solvability of the equations of motion of multicomponent viscous compressible fluids (homogeneous mixtures of fluids, multifluids). Concerning the origin of the model and its physical interpretation, we refer the reader to \cite{mamprok.jfcs2}. An overview of the options for formulating the model and the known results can be found in \cite{mamprok.france} and \cite{mamprok.semi17}. Related multi-velocity models of multifluids are considered in \cite{mamprok.dorovski}, \cite{mamprok.giovang}, \cite{mamprok.nigm} and \cite{mamprok.raj}. As the first results on the well-posedness of the multidimensional equations of multifluids, we can refer to \cite{mamprok.frehsegm1}, \cite{mamprok.frehsegm2} and \cite{mamprok.frehsew}. Solvability for related models is shown in \cite{mamprok.jms}, \cite{mamprok.muchapz}, \cite{mamprok.semiprok17} and \cite{mamprok.jfcsprok}.

Weak solutions for multidimensional barotropic problems for the model con\-si\-dered in the paper are constructed in the steady version in \cite{mamprok.semi162} and \cite{mamprok.jaimprok} (polytropic case), and then in \cite{mamprok.smz1} and \cite{mamprok.smz2} (general case); in the unsteady version in \cite{mamprok.semi161} (polytropic case), and then in \cite{mamprok.izvran18} (general case). Similar results for the heat-conductive model are obtained in \cite{mamprok.izvran14}. For a number of reasons, including the purpose of constructing more regular solutions, one-dimensional formulations are of interest. It should be noted that the dimension with respect to the number of components (constituents) of a multifluid is not logically and technically related to the number of the spatial variables, and the interaction of the constituents via the viscous terms transforms the system of differential equations governing the motion of a multifluid into a system essentially different from the one-component system. Therefore, despite the developed theory of one-dimensional flows of a viscous gas (see \cite{mamprok.kazhshel} for example), the one-dimensional theory of multifluids remains relevant.

The specificity of the paper is that we consider a variant of the model with an average velocity in the transport operator.

\section{Statement of the problem}

\noindent\indent
Let us consider the system of equations governing motions of multicomponent viscous compressible fluids without taking into account chemical reactions in the case of one spatial variable:
$$\partial_{t}\rho_{i}+\partial_{x}(\rho_{i} v)=0,$$
$$\rho_{i}\left(\partial_{t}u_{i}+v\partial_{x}u_{i}\right)=\partial_{x}P_{i},\quad i=1,\ldots,N.$$
Here $N\geqslant 2$ is the number of components, $\rho_{i}$ is the density if the $i$-th constituent, $u_{i}$ is the velocity of the $i$-th component,
$\displaystyle v=\frac{1}{N}\sum\limits_{i=1}^{N}u_{i}$ is the average velocity of the multifluid, and $P_{i}$ are the stresses. We accept generalized Newton's hypothesis
$$P_{i}=-p+\sum\limits_{j=1}^{N}\mu_{ij}\partial_{x}u_{j},$$
where $p$ is the pressure, and the viscosity coefficients $\{\mu_{ij}\}_{i, j=1}^{N}$ form the symmetric matrix~$\textbf{M}$. Moreover, $\textbf{M}>0$, i.~e.
$(\textbf{M}\boldsymbol{\xi},\boldsymbol{\xi})\geqslant C_{0}(\textbf{M})|\boldsymbol{\xi}|^2$ for all $\boldsymbol{\xi}\in{\mathbb R}^N$ with a constant $C_{0}(\textbf{M})>0$.

The written equations together with the constitutive equation
$$p=K\rho^{\gamma},\quad \rho=\sum\limits_{i=1}^{N}\rho_{i},\quad K={\rm const}>0,\quad \gamma={\rm const}>1$$
form a closed system
\begin{align}\label{newcontinuity.double}\partial_{t}\rho_{i}+\partial_{x}(\rho_{i} v)=0,\quad v=\frac{1}{N}\sum\limits_{i=1}^{N}u_{i},\end{align}
\begin{align}\label{newmomentum.double}\rho_{i}\left(\partial_{t}u_{i}+v\partial_{x}u_{i}\right)+K\partial_{x} \rho^{\gamma}=
\sum\limits_{j=1}^N \mu_{ij}\partial_{xx}u_{j},\quad i=1,\ldots,N,\quad \rho=\sum\limits_{i=1}^{N}\rho_{i}.\end{align}
Let us consider this system in the rectangular $Q_{T}$ (here and below\linebreak $Q_{t}=(0, 1)\times(0, t)$) with an arbitrary finite height $T$, $0<T<\infty$, and endow this system with the following initial and boundary conditions ($i=1,\ldots,N$):
\begin{align}\label{nachusl.double}\rho_{i}|_{t=0}=\rho_{0i}(x), \quad u_{i}|_{t=0}=u_{0i}(x),\quad x\in [0, 1],\end{align}
\begin{align}\label{boundvelocity.double}u_{i}|_{x=0}=u_{i}|_{x=1}=0,\quad t\in [0, T].\end{align}

{\bfseries Definition 1.} By a strong solution to the problem \eqref{newcontinuity.double}--\eqref{boundvelocity.double} we mean a collection of $2N$ functions
$(\rho_{1},\ldots, \rho_{N}, u_{1},\ldots, u_{N})$ such that the equations \eqref{newcontinuity.double} and \eqref{newmomentum.double} are satisfied a.~e. in $Q_{T}$, the initial conditions \eqref{nachusl.double} are satisfied for a.~a. $x\in (0, 1)$, the boundary conditions \eqref{boundvelocity.double} are satisfied for a.~a. $ t\in (0, T)$, and the following inequalities and inclusions hold ($i=1,\ldots,N$)
$$\rho_{i}>0,\quad \rho_{i}\in L_{\infty}\big(0, T; W^{1}_{2}(0, 1)\big), \quad  \partial_{t}\rho_{i}\in L_{\infty}\big(0, T; L_{2}(0, 1)\big),$$
$$u_{i}\in L_{\infty}\big(0, T; W^{1}_{2}(0, 1)\big)\bigcap L_{2}\big(0, T; W^{2}_{2}(0, 1)\big),\quad \partial_{t}u_{i} \in L_{2}(Q_{T}).$$

\section{Main result}

\noindent\indent
{\bfseries Theorem 2.} Suppose that the initial data in \eqref{nachusl.double} satisfy the conditions
$$\rho_{0i}\in W^{1}_{2}(0,1),\quad \rho_{0i}>0,\quad u_{0i}\in \overset{\circ}{W^1_2}(0, 1),\quad i=1,\ldots,N\quad (N\geqslant 2),$$
the symmetric viscosity matrix $\textbf{M}$ is positive definite, and the polytropic index $\gamma>1$ as well as the constants $0<K,T<\infty$ are given. Then there exists a unique strong solution to the problem \eqref{newcontinuity.double}--\eqref{boundvelocity.double} in the sense of Definition~1.

Since the uniqueness and the corresponding local result are obtained in \cite{mamprok.local}, the proof of Theorem~2 reduces to obtaining global a priori estimates, which is the main content of the paper.

\section{Lagrangian coordinates}

\noindent\indent
During the study of the problem \eqref{newcontinuity.double}--\eqref{boundvelocity.double}, the parallel use of the Lagrangian coordinates is convenient. Let us accept
$\displaystyle y(x,t)=\int\limits_{0}^{x}\rho(s,t)\,ds$ and~$t$ as new independent variables. Then the system \eqref{newcontinuity.double}, \eqref{newmomentum.double} takes the form
\begin{align}\label{1newcontinuity1lagr}
\partial_{t}\rho_{i}+\rho\rho_{i}\partial_{y}v=0,\quad v=\frac{1}{N}\sum\limits_{i=1}^{N}u_{i},
\end{align}
\begin{align}\label{1newmomentum1lagr}
\frac{\rho_{i}}{\rho}\partial_{t}u_{i}+K\partial_{y} \rho^{\gamma}=
\sum\limits_{j=1}^N \mu_{ij}\partial_{y}(\rho\partial_{y}u_{j}),\quad i=1,\ldots,N,\quad \rho=\sum\limits_{i=1}^{N}\rho_{i},
\end{align}
the domain $Q_{T}$ is transformed into the rectangular $\Pi_{T}=(0, d)\times(0, T)$, where $\displaystyle d=\int\limits_{0}^{1}\rho_{0}\,dx>0$, $\displaystyle \rho_{0}=\sum\limits_{i=1}^{N}\rho_{0i}$, and the initial and boundary conditions take the form ($i=1,\ldots,N$)
\begin{align}\label{1nachusl1lagr}
\rho_{i}|_{t=0}=\widetilde{\rho}_{0i}(y), \quad u_{i}|_{t=0}=\widetilde{u}_{0i}(y),\quad y\in [0, d],
\end{align}
\begin{align}\label{1boundvelocity1lagr}
u_{i}|_{y=0}=u_{i}|_{y=d}=0,\quad t\in [0, T].
\end{align}

\section{Estimates of the concentrations}

\noindent\indent
Let us consider a hypothetical solution $(\rho_{1},\ldots, \rho_{N}, u_{1},\ldots, u_{N})$ to the problem \eqref{newcontinuity.double}--\eqref{boundvelocity.double} which possesses all necessary differential properties, and such that the densities $\rho_{i}$, $i=1,\ldots,N$, are positive and bounded (see Definition~1).

First of all, we note that the summation of \eqref{1newcontinuity1lagr} with respect to $i=1,\ldots,N$ gives
\begin{align}\label{1602171}
\partial_{t}\rho+\rho^{2}\partial_{y}v=0,
\end{align}
and hence
$$\partial_{t}\left(\frac{\rho_{i}}{\rho}\right)=0, \quad i=1,\ldots,N.$$
Hence, due to \eqref{1nachusl1lagr} we get the equalities
\begin{align}\label{1602173}
\frac{\rho_{i}(y,t)}{\rho(y,t)}=\frac{\widetilde{\rho}_{0i}(y)}{\widetilde{\rho}_{0}(y)}\quad\text{as}\quad(y,t)\in[0, d]\times[0, T]
\end{align}
for all $i=1,\ldots,N$, where $\displaystyle \widetilde{\rho}_{0}=\sum\limits_{i=1}^{N}\widetilde{\rho}_{0i}$. In the Eulerian coordinates $(x,t)$, the ratios
$\displaystyle \frac{\rho_{i}}{\rho}$ satisfy the transport equations, and we only have the inequalities
\begin{multline}\label{2602171}
\inf\limits_{[0,1]}\frac{\rho_{0i}(x)}{\rho_{0}(x)}\leqslant\frac{\rho_{i}(x,t)}{\rho(x,t)}\leqslant\sup\limits_{[0,1]}\frac{\rho_{0i}(x)}{\rho_{0}(x)}\leqslant 1\quad\text{as}\quad(x,t)\in[0, 1]\times[0, T],\\ i=1,\ldots,N.
\end{multline}

\section{First a priori estimates}

\noindent\indent
Typically for the compressible Navier---Stokes theory, the energy inequality immediately entails the estimates for the kinetic energy, the rate of energy dissipation and the potential energy of the multifluid constituents, as we show below.

{\bfseries Lemma 3.}  Under the assumptions of Theorem 2, there exists a positive constant\footnote{Hereinafter, $C$ with indices denotes positive constants which depend on the initial data, physical constants and $T$.}
$\displaystyle C_{1}\left(\left\{\inf\limits_{[0,1]}\frac{\rho_{0i}}{\rho_{0}}\right\}, \left\{\|\sqrt{\rho_{0i}}u_{0i}\|_{L_{2}(0,1)}\right\},
\|\rho_{0}\|_{L_{\gamma}(0,1)}, K, \textbf{M}, N, \gamma\right)$, such that the following estimate holds
\begin{equation}\label{lemma1}
\sum\limits_{i=1}^{N}\left(\|\sqrt{\rho}u_{i}\|_{L_{\infty}\big(0, T;L_{2}(0, 1)\big)}
+\|\partial_{x}u_{i}\|_{L_{2}(Q_{T})}\right)+\|\rho\|_{L_{\infty}\big(0, T;L_{\gamma}(0,1)\big)}\leqslant C_{1}.
\end{equation}

{\bfseries Proof.} Let us multiply the equations \eqref{newmomentum.double} by $u_{i}$, integrate over $(0, 1)$ and sum with respect to $i=1,\ldots,N$. In view of \eqref{newcontinuity.double}, \eqref{boundvelocity.double} and $\textbf{M}>0$, the following relations hold
$$\sum\limits_{i=1}^{N}\int\limits^{1}_{0}\Big(\rho_{i}\partial_{t}u_{i}+\rho_{i} v\partial_{x}u_{i}\Big)u_{i}\,dx= \frac{1}{2}\frac{d}{dt}\sum\limits_{i=1}^{N}
\int\limits^{1}_{0}\rho_{i} u_{i}^{2}\, dx,$$
$$\sum\limits_{i=1}^{N}\int\limits^{1}_{0}u_{i}K\partial_{x}\rho^{\gamma}\,dx=-KN\int\limits^{1}_{0}\rho^{\gamma}\partial_{x}v\, dx=
\frac{KN}{\gamma-1}\frac{d}{dt}\int\limits^{1}_{0}\rho^{\gamma}\, dx,$$
\begin{multline}\label{lemma51.4}
\sum\limits_{i,j=1}^N \mu_{ij}\int\limits\limits_{0}^{1}(\partial_{xx}u_{j})u_{i}\,
dx=-\sum\limits_{i,j=1}^N \mu_{ij}\int\limits^{1}_{0}(\partial_{x}u_{i})(\partial_{x}u_{j})\,dx\leqslant\\ \leqslant -C_{0}(\textbf{M})
\sum\limits_{i=1}^N\int\limits\limits_{0}^{1}|\partial_{x}u_{i}|^{2}\,dx,
\end{multline}
and hence we get the inequality
\begin{equation}\label{lemma51.6}
\frac{d}{dt}\sum\limits_{i=1}^{N}\int\limits^{1}_{0}\left(\frac{1}{2}\rho_{i} u_{i}^{2}+\frac{K}{\gamma-1}\rho^{\gamma}\right)\,dx+
C_{0}\sum\limits_{i=1}^N\int\limits\limits_{0}^{1}|\partial_{x}u_{i}|^{2}\,dx\leqslant 0.
\end{equation}
We integrate \eqref{lemma51.6} over $(0, t)$, and using \eqref{nachusl.double}, we obtain the bound
$$\sum\limits_{i=1}^{N}\int\limits^{1}_{0}\left(\frac{1}{2}\rho_{i} u_{i}^{2}+\frac{K}{\gamma-1}\rho^{\gamma}\right)\,dx+
C_{0}\sum\limits_{i=1}^N\int\limits\limits_{0}^{t}\int\limits\limits_{0}^{1}|\partial_{x}u_{i}|^{2}\,dxd\tau\leqslant$$
$$\leqslant\sum\limits_{i=1}^{N}\int\limits^{1}_{0}\left(\frac{1}{2}\rho_{0i}u_{0i}^{2}+\frac{K}{\gamma-1}\rho_{0}^{\gamma}\right)\,dx,$$
which, due to \eqref{2602171}, implies the conclusion of Lemma 3.

{\bfseries Remark 4.} In the Lagrangian coordinates, the bound \eqref{lemma1} takes the form
\begin{multline}\label{lemma1lagr}
\sum\limits_{i=1}^{N}\left(\left\|u_{i}\right\|_{L_{\infty}\big(0, T;L_{2}(0, d)\big)}+\|\sqrt{\rho}\partial_{y}u_{i}\|_{L_{2}(\Pi_{T})}\right)+
\|\rho\|_{L_{\infty}\big(0, T;L_{\gamma-1}(0,d)\big)}\leqslant\\ \leqslant C_{2}(C_{1},\gamma).
\end{multline}

{\bfseries Remark 5.} In view of \eqref{1boundvelocity1lagr}, we obviously get the following inequality from~\eqref{lemma1}
\begin{equation}\label{lemma1lagr2602171}
\sum\limits_{i=1}^{N}\left\|u_{i}\right\|_{L_{2}\big(0, T;L_{\infty}(0, 1)\big)}\leqslant C_{1}.
\end{equation}

\section{Bound of the density from above}

\noindent\indent
The crucial a priori estimates are the bounds of strict positiveness and boundedness of the densities of the multifluid constituents. First we prove the bound for the densities from above.

{\bfseries Lemma 6.} There exists a constant\\
$C_{3}\left(C_{1},  \|\rho_{0}\|_{L_{\infty}(0,1)}, \{\|\rho_{0i}u_{0i}\|_{L_{1}(0, 1)}\}, K, \textbf{M}, N, T, d, \gamma\right)$ such that
\begin{equation}\label{supernew1602171}
\rho(x,t)\leqslant C_{3}\quad\text{as}\quad(x,t)\in[0, 1]\times[0, T].
\end{equation}

{\bfseries Proof.} Let us rewrite the equations \eqref{newmomentum.double}, using \eqref{newcontinuity.double}, in the form
\begin{multline}\label{new1602171}
\frac{1}{N}\sum\limits_{j=1}^N \widetilde{\mu}_{ij}\left(\partial_{t}(\rho_{j}u_{j})+\partial_{x}(\rho_{j}v u_{j})\right)+\frac{K}{N}\left(\sum\limits_{j=1}^N \widetilde{\mu}_{ij}\right)\partial_{x} \rho^{\gamma}=\frac{1}{N}\partial_{xx}u_{i},\\ i=1,\ldots,N,
\end{multline}
where $\widetilde{\mu}_{ij}$ are the entries of the matrix $\widetilde{\textbf{M}}=\textbf{M}^{-1}>0$, and then sum \eqref{new1602171} with respect to $i=1,\ldots,N$, then we get
\begin{align}\label{new1602172}
\partial_{t}V=\partial_{x}\left(\partial_{x}v-\widetilde{K}\rho^{\gamma}-vV\right),
\end{align}
where $\displaystyle V=\frac{1}{N}\sum\limits_{i,j=1}^N \widetilde{\mu}_{ij}\rho_{j}u_{j}$, $\displaystyle \widetilde{K}=\frac{K}{N}\sum\limits_{i, j=1}^{N}\widetilde{\mu}_{ij}>0$.
We denote
\begin{equation}\label{lemma52.1}
\alpha(x,t)=\int\limits_{0}^t\Big(\partial_{x}v-\widetilde{K}\rho^{\gamma}- v V\Big)\,d\tau+\int\limits_{0}^xV_{0}\,ds,
\end{equation}
where $V_{0}(x)=V(x,0)$. In view of \eqref{lemma1}, we have
$$\|\partial_{x}\alpha\|_{L_{\infty}\big(0, T;L_{1}(0, 1)\big)}=\|V\|_{L_{\infty}\big(0, T;L_{1}(0, 1)\big)}\leqslant C_{4}(C_{1}, \textbf{M}, d),$$
$$\sup\limits_{[0, T]}\left|\int\limits^{1}_{0}\alpha\,dx\right| \leqslant T\max\limits_{1\leqslant i,j\leqslant N}|\widetilde{\mu}_{ij}|\sum\limits_{i=1}^{N}\sup\limits_{[0, T]}\int\limits^{1}_{0}\rho u_{i}^{2}\, dx+\widetilde{K}T\sup\limits_{[0, T]}\int\limits^{1}_{0}\rho^{\gamma}\, dx+$$
$$+\max\limits_{1\leqslant i, j\leqslant N}|\widetilde{\mu}_{ij}| \sum\limits_{i=1}^{N}\int\limits_{0}^1\rho_{0i}|u_{0i}|\, dx\leqslant
C_{5}\left(C_{1}, \{\|\rho_{0i}u_{0i}\|_{L_{1}(0, 1)}\},\textbf{M}, \widetilde{K}, T, \gamma\right),$$
and hence, using Poincar\'e's inequality, we get\footnote{The bound $C_{6}$ depends on the size of the flow domain.}
$$\sup\limits_{[0, T]}\int\limits^{1}_{0}|\alpha|\, dx \leqslant \sup\limits_{[0, T]}\int\limits^{1}_{0}|\partial_{x}\alpha |\,dx+
\sup\limits_{[0, T]}\left|\int\limits^{1}_{0}\alpha \,dx\right|\leqslant C_{6}(C_{4}, C_{5}),$$
and we arrive at the boundedness of $\alpha$ in $L_{\infty}\big(0,T;W_{1}^{1}(0, 1)\big)$. Using this and the fact $W_{1}^{1}(0, 1)\hookrightarrow L_{\infty}(0,1)$, we obtain the estimate
$$\|\alpha\|_{L_{\infty}(Q_{T})}\leqslant C_{7}\left(C_{4}, C_{6}\right).$$
Let us note that, in view of \eqref{newcontinuity.double}, \eqref{new1602172} and \eqref{lemma52.1}, the following relations hold
$$d_{t}(\rho e^{\alpha})=-\widetilde{K}e^{\alpha}\rho^{\gamma+1}\leqslant 0,\quad \text{where}\quad d_{t}=\partial_{t}+v\partial_{x},$$
and hence
$$\rho e^{\alpha}\leqslant\sup\limits_{[0,1]}{\rho_{0}}\exp\left(\int\limits_{0}^{1}|V_{0}|\,dx\right),$$
so that we arrive at the conclusion of Lemma 6.

\section{The bound for the derivative of the density}

\noindent\indent
In order to obtain the bound for the densities from below, we first need to prove the boundedness of the first spatial derivative of the logarithm of the total density. Specifically, the following assertion holds.

{\bfseries Lemma 7.} There exists a constant
$$\displaystyle C_{8}\left(C_{1}, C_{2}, C_{3}, \left\{\|\widetilde{u}_{0i}\|_{L_{2}(0,d)}\right\}, \left\{\left\|\frac{\widetilde{\rho}_{0i}}{\widetilde{\rho}_{0}}\right\|_{W^{1}_{2}(0,d)}\right\},
\|\left(\ln{\widetilde{\rho}}_{0}\right)^{\prime}\|_{L_{2}(0,d)}, \textbf{M}, N\right)$$
such that
\begin{align}\label{eq01021710}
\|\partial_{y}\ln\rho\|_{L_{\infty}\big(0, T; L_{2}(0, d)\big)}\leqslant C_{8}.
\end{align}

{\bfseries Proof.} Let us use the equations in the form \eqref{1newcontinuity1lagr}, \eqref{1newmomentum1lagr}. We rewrite the equations \eqref{1newmomentum1lagr} as
\begin{align}\label{newnewmomentum.double}
\frac{1}{N}\sum\limits_{j=1}^{N}\widetilde{\mu}_{ij}\frac{\rho_{j}}{\rho}\partial_{t}u_{j}+\frac{K}{N}\left(\sum\limits_{j=1}^{N}\widetilde{\mu}_{ij}\right)\partial_{y} \rho^{\gamma}=
\frac{1}{N}\partial_{y}(\rho\partial_{y}u_{i}),\quad i=1,\ldots,N,
\end{align}
and then sum \eqref{newnewmomentum.double} with respect to $i=1,\ldots,N$, then we get, using \eqref{1602173}, the resulting relation
\begin{align}\label{newnewmomentum1newnew}
\frac{1}{N}\sum\limits_{i, j=1}^{N}\widetilde{\mu}_{ij}\frac{\widetilde{\rho}_{0j}}{\widetilde{\rho}_{0}}\partial_{t}u_{j}+\widetilde{K}\partial_{y}\rho^{\gamma}=
\partial_{y}\left(\rho\partial_{y}v\right).
\end{align}
We extract from \eqref{1602171} that
\begin{align}\label{eq0102173}
\rho\partial_{y}v=-\partial_{t}\ln\rho
\end{align}
and substitute this into \eqref{newnewmomentum1newnew}, then we get
$$\partial_{ty}\ln\rho+\widetilde{K}\partial_{y}\rho^{\gamma}=-\frac{1}{N}\sum\limits_{i,j=1}^{N}\widetilde{\mu}_{ij}
\frac{\widetilde{\rho}_{0j}}{\widetilde{\rho}_{0}}\partial_{t}u_{j}.$$
We multiply this equality by $\displaystyle \partial_{y}\ln\rho=:w$ and integrate over $y\in (0,d)$, then we obtain
\begin{align}\label{eq0102172}
\frac{1}{2}\frac{d}{dt}\left(\int\limits_{0}^{d}w^{2}\, dy\right)+\widetilde{K}\gamma\int\limits_{0}^{d}\rho^{\gamma}w^{2}\,dy
=-\frac{1}{N}\sum\limits_{i, j=1}^{N}\widetilde{\mu}_{ij}\int\limits_{0}^{d}\left(\frac{\widetilde{\rho}_{0j}}{\widetilde{\rho}_{0}}\partial_{t}u_{j}\right)w\,dy.
\end{align}
Let us transform the right-hand side of \eqref{eq0102172} via the integration by parts and using \eqref{eq0102173}:
\begin{align}\label{eq01021777}\begin{array}{c}\displaystyle
-\frac{1}{N}\sum\limits_{i, j=1}^{N}\widetilde{\mu}_{ij}\int\limits_{0}^{d}\left(\frac{\widetilde{\rho}_{0j}}{\widetilde{\rho}_{0}}\partial_{t}u_{j}\right)w\,dy
=-\frac{d}{dt}\left(\frac{1}{N}\sum\limits_{i, j=1}^{N}\widetilde{\mu}_{ij}\int\limits_{0}^{d} \frac{\widetilde{\rho}_{0j}}{\widetilde{\rho}_{0}}u_{j}w\,dy\right)+\\ \\
\displaystyle
+\frac{1}{N}\sum\limits_{i, j=1}^{N}\widetilde{\mu}_{ij}\int\limits_{0}^{d} \rho u_{j}(\partial_{y}v)\left(\frac{\widetilde{\rho}_{0j}}{\widetilde{\rho}_{0}}\right)^{\prime}\,dy+
\frac{1}{N}\sum\limits_{i, j=1}^{N}\widetilde{\mu}_{ij}\int\limits_{0}^{d}\frac{\widetilde{\rho}_{0j}\rho}{\widetilde{\rho}_{0}}(\partial_{y}v)(\partial_{y}u_{j})\,dy.
\end{array}
\end{align}
Thus, after integration of \eqref{eq0102172} with respect to $t$, taking into account \eqref{supernew1602171} and \eqref{eq01021777}, we get
$$\|w\|^{2}_{L_{2}(0, d)}+2\widetilde{K}\gamma\int\limits_{0}^{t}\int\limits_{0}^{d}\rho^{\gamma}w^{2}\,dyd\tau\leqslant$$
$$\leqslant\|w_{0}\|^{2}_{L_{2}(0, d)}-\frac{2}{N}\sum\limits_{i, j=1}^{N}\widetilde{\mu}_{ij}\int\limits_{0}^{d} \frac{\widetilde{\rho}_{0j}}{\widetilde{\rho}_{0}}u_{j}w\,dy+
\frac{2}{N}\sum\limits_{i, j=1}^{N}\widetilde{\mu}_{ij}\int\limits_{0}^{d} \frac{\widetilde{\rho}_{0j}}{\widetilde{\rho}_{0}}\widetilde{u}_{0j}w_{0}\,dy+$$
$$+\frac{2\sqrt{C_{3}}}{N}\sum\limits_{i, j=1}^{N}|\widetilde{\mu}_{ij}|\int\limits_{0}^{t}\left\|\left(\frac{\widetilde{\rho}_{0j}}{\widetilde{\rho}_{0}}\right)^{\prime}\right\|_{L_{2}(0, d)}\|u_{j}\|_{L_{\infty}(0, d)}\|\sqrt{\rho}\partial_{y}v\|_{L_{2}(0, d)}\,d\tau+$$
$$+\frac{2}{N}\sum\limits_{i, j=1}^{N}|\widetilde{\mu}_{ij}|\sup\limits_{[0,d]}\frac{\widetilde{\rho}_{0j}}{\widetilde{\rho}_{0}}\int\limits_{0}^{t}\|\rho (\partial_{y}v)(\partial_{y}u_{j})\|_{L_{1}(0, d)}\,d\tau,$$
where $w_{0}=(\ln{\widetilde{\rho}_{0}})^{\prime}$. Using the estimates \eqref{lemma1lagr} and \eqref{lemma1lagr2602171}, we derive from this the inequality~\eqref{eq01021710}, concluding the proof of Lemma~7.

\section{The estimate of the density from below}

\noindent\indent
In this section, we finish obtaining the crucial estimates of strict positiveness and boundedness of the densities via the following assertion.

{\bfseries Lemma 8.} There exists a constant $C_{9}(C_{8},d)$ such that
\begin{align}\label{eq01021714}
\rho(y,t)\geqslant C_{9}\quad\text{as}\quad(y,t)\in[0, d]\times[0, T].
\end{align}

{\bfseries Proof.} The continuity equation for $\rho$ immediately leads, for any\linebreak $t\in[0,T]$, to the existence of a point $z(t)\in[0, d]$ such that
\begin{align}\label{eq0102176}
\rho(z(t),t)=d.
\end{align}
Hence, we can use the representation
$$\ln\rho(y,t)=\ln\rho(z(t),t)+\int\limits_{z(t)}^{y}\partial_{s}\ln\rho(s,t)\, ds,$$
from which, via H\"older's inequality, and using \eqref{eq01021710} and \eqref{eq0102176}, we get
$$|\ln\rho(y,t)|\leqslant |\ln{d}|+\sqrt{d}\|w\|_{L_{2}(0,d)}\leqslant C_{10}(C_{8},d).$$
This leads immediately to \eqref{eq01021714}, and Lemma 8 is proved.

{\bfseries Remark 9.} The equalities \eqref{1602173} and the estimates in Lemmas 6 and 8 imply that for all $i=1,\ldots,N$ we have
\begin{align}\label{eq01021714111}
C_{11}\leqslant\rho_{i}(y,t)\leqslant C_{3}\quad\text{as}\quad(y,t)\in[0, d]\times[0, T],
\end{align}
where $\displaystyle C_{11}=C_{11}\left(C_{9},\left\{\inf\limits_{[0,d]}\frac{\widetilde{\rho}_{0i}(y)}{\widetilde{\rho}_{0}(y)}\right\}\right)$.

\section{Further bounds}

\noindent\indent
Concluding the proof of Theorem 2, we obtain the estimates for the derivatives of the densities and velocities of the multifluid constituents.

{\bfseries Remark 10.} From the estimates in Lemmas 6 and 7 and the formula~\eqref{1602173}, it follows that
\begin{align}\label{finalocenla}
\left\|\partial_{x}\rho_{i}\right\|_{L_{\infty}(0, T; L_{2}(0,1))}\leqslant C_{12},\quad i=1,\ldots,N,
\end{align}
where $\displaystyle C_{12}=C_{12}\left(C_{3}, C_{8}, \left\{\left\|\left(\frac{\widetilde{\rho}_{0}}{\widetilde{\rho}_{0i}}\right)^{\prime}\right\|_{L_{2}(0,d)}\right\}, \left\{\sup\limits_{[0,d]}\frac{\widetilde{\rho}_{0i}}{\widetilde{\rho}_{0}}\right\}\right)$.

\medskip

{\bfseries Lemma 11.} There exists\\ $C_{13}\left(C_{1}, C_{3}, C_{11}, C_{12}, \{\|u_{0i}^{\prime}\|_{L_{2}(0, 1)}\}, K, \textbf{M}, N, T, \gamma\right)$ such that
$$\sum\limits_{i=1}^{N}\left(\|\partial_{x}u_{i}\|_{L_{\infty}(0, T; L_{2}(0,1))}+\|\partial_{xx}u_{i}\|_{L_{2}(Q_{T})}+\|\partial_{t}u_{i}\|_{L_{2}(Q_{T})}\right)\leqslant C_{13}.$$

{\bfseries Proof.} We derive from \eqref{finalocenla} that
\begin{align}\label{eq01021715}
\left\|(\partial_{x}\rho)(t)\right\|_{L_{2}(0,1)}\leqslant C_{14}(C_{12}, N)\quad \forall\, t\in[0,T].
\end{align}
Using the idea of \cite{mamprok.basovshel}, we square the momentum equations \eqref{newmomentum.double} and sum the result with respect to $i=1,\ldots,N$, then we get
\begin{multline}\label{1303171}\sum\limits_{i=1}^{N}\rho_{i}(\partial_{t}u_{i})^{2}+\sum\limits_{i=1}^{N}\frac{1}{\rho_{i}}
\left(\sum\limits_{j=1}^{N}\mu_{ij}\partial_{xx}u_{j}\right)^{2}-2\sum\limits_{i=1}^{N}(\partial_{t}u_{i})\left(\sum\limits_{j=1}^{N}\mu_{ij}\partial_{xx}u_{j}\right)=\\
=\sum\limits_{i=1}^{N}\rho_{i}\left(\frac{K\partial_{x}\rho^{\gamma}}{\rho_{i}}+v\partial_{x}u_{i}\right)^{2}.\end{multline}
Let us introduce a function $\beta(t)$ via the relation
$$\beta(t)=\sum\limits_{i, j=1}^{N}\mu_{ij}\int\limits_{0}^{1}(\partial_{x}u_{i})(\partial_{x}u_{j})\, dx+$$
$$+\sum\limits_{i=1}^{N}\int\limits_{0}^{t}\int\limits_{0}^{1}\left(\rho_{i}(\partial_{t}u_{i})^{2}+
\frac{1}{\rho_{i}}\left(\sum\limits_{j=1}^{N}\mu_{ij}\partial_{xx}u_{j}\right)^{2}\right)\, dxd\tau.$$
Then \eqref{1303171} and the inequalities \eqref{lemma51.4}, \eqref{supernew1602171}, \eqref{eq01021714111} and \eqref{eq01021715} give the
estimate\footnote{Here the symmetry of the matrix $\textbf{M}$ is used.}
$$\beta^{\prime}(t)\leqslant C_{15}+C_{16}\left(\sum\limits_{j=1}^{N}\|u_{j}\|^{2}_{L_{\infty}(0,1)}\right)\left(\sum\limits_{i, j=1}^{N}\mu_{ij}
\int\limits_{0}^{1}(\partial_{x}u_{i})(\partial_{x}u_{j})\, dx\right)\leqslant$$
$$\leqslant C_{15}+C_{16}\left(\sum\limits_{j=1}^{N}\|u_{j}\|^{2}_{L_{\infty}(0,1)}\right)\beta(t),$$
where $C_{15}=C_{15}(C_{3}, C_{11}, C_{14}, K, N, \gamma)$, $C_{16}=C_{16}(C_{3}, \textbf{M})$, from which, via Gronwall's lemma (see also \eqref{lemma1lagr2602171}), it follows that
$$\beta(t)\leqslant C_{17}\left(C_{1}, C_{15}, C_{16}, \{\|u_{0i}^{\prime}\|_{L_{2}(0, 1)}\}, \textbf{M}, T\right),$$
and we arrive at the conclusion of Lemma 11.

{\bfseries Remark 12.} It follows immediately from the continuity equations \eqref{newcontinuity.double} and the estimates in Lemmas 6 and 11 and Remark 10, that
$$\|\partial_{t}\rho_{i}\|_{L_{\infty}(0, T; L_{2}(0,1))}\leqslant C_{18}(C_{3}, C_{12}, C_{13}),\quad i=1,\ldots,N.$$

\section{Proof of Theorem 2}

\noindent\indent
Basing on the global a priori estimates proved in Sections~5--10, we can continue the local solution (obtained in Theorem 2 from \cite{mamprok.local}) into the entire~$Q_{T}$ (see, e.~g., \cite{mamprok.akm}, P.~40, or \cite{mamprok.weigant92}, P.~20). The uniqueness of this solution is shown in \cite{mamprok.local}. Thus, Theorem~2 is proved.


\newpage
\begin{thebibliography}{99}
\bibitem{mamprok.akm} S.~N.~Antontsev, A.~V.~Kazhikhov and V.~N.~Monakhov, {\it Boundary value problems in mechanics of nonhomogeneous fluids},
Studies in Mathematics and its Applications, {\bfseries 22}, North-Holland Publishing Co., Amsterdam, 1990.

\bibitem{mamprok.basovshel} I.~V.~Basov and V.~V.~Shelukhin, {\it On equations of a nonlinear compressible fluid with a discontinuous constitutive law}, Siberian Math. J.,
{\bfseries 40}:3 (1999), 435--445.

\bibitem{mamprok.dorovski} V.~N.~Dorovsky and Yu.~V.~Perepechko, {\it Theory of the partial melting}, Sov. Geology and Geophysics, 9 (1989), 56--64 (in Russian).

\bibitem{mamprok.frehsegm1} J.~Frehse, S.~Goj and J.~Malek, {\it On a Stokes-like system for mixtures of fluids}, SIAM J. Math. Anal., {\bfseries 36}:4 (2005), 1259--1281.

\bibitem{mamprok.frehsegm2} J.~Frehse, S.~Goj and J.~Malek, {\it A uniqueness result for a model for mixtures in the absence of external forces and interaction momentum},
Appl. Math., {\bfseries 50}:6 (2005), 527--541.

\bibitem{mamprok.frehsew} J.~Frehse and W.~Weigant, {\it On quasi-stationary models of mixtures of compressible fluids}, Appl. Math., {\bfseries 53}:4 (2008), 319--345.

\bibitem{mamprok.giovang} V.~Giovangigli, {\it Multicomponent flow modeling}, Birkh\"auser, Boston, 1999.

\bibitem{mamprok.kazhshel} A.~V.~Kazhikhov and V.~V.~Shelukhin, {\it Unique global solution with respect to time of initial-boundary value problems for one-dimensional equations of a viscous gas}, Journal of Applied Mathematics and Mechanics, {\bfseries 41}:2 (1977), 273--282.

\bibitem{mamprok.france} A.~E.~Mamontov and D.~A.~Prokudin, {\it Viscous compressible multi-fluids: modeling and multi-D existence}, Methods and applications of analysis,
{\bfseries 20}:2 (2013), 179--195.

\bibitem{mamprok.izvran14} A.~E.~Mamontov and D.~A.~Prokudin, {\it Solubility of a stationary boundary-value problem for the equations of motion of a one-temperature mixture of viscous compressible heat--conducting fluids}, Izvestiya: Mathematics, {\bfseries 78}:3 (2014), 554--579.

\bibitem{mamprok.smz1} A.~E. Mamontov and D.~A. Prokudin, {\it Solvability of the regularized steady problem of the spatial motions of multicomponent viscous compressible fluids},
Siberian Math. J., {\bfseries 57}:6 (2016), 1044--1054.

\bibitem{mamprok.semi161} A.~E. Mamontov and D.~A. Prokudin, {\it Solubility of initial boundary value problem for the equations of polytropic motion of multicomponent viscous compressible fluids}, Siberian Electr. Math. Reports, {\bfseries 13} (2016), 541--583 (in Russian).

\bibitem{mamprok.semi162} A.~E. Mamontov and D.~A. Prokudin, {\it Solubility of steady boundary value problem for the equations of polytropic motion of multicomponent viscous compressible fluids}, Siberian Electr. Math. Reports, {\bfseries 13} (2016), 664--693 (in Russian).

\bibitem{mamprok.smz2} A.~E. Mamontov and D.~A. Prokudin, {\it Existence of weak solutions to the three-dimensional problem of steady barotropic motions of mixtures of viscous compressible fluids}, Siberian Math. J., {\bfseries 58}:1 (2017), 113--127.

\bibitem{mamprok.semi17} A.~E.~Mamontov and D.~A.~Prokudin, {\it Viscous compressible homogeneous multi-fluids with multiple velocities: barotropic existence theory}, Siberian Electr. Math. Reports, {\bfseries 14} (2017), 388--397.

\bibitem{mamprok.jfcs2} A.~E.~Mamontov and D.~A.~Prokudin, {\it Modeling viscous compressible barotropic multi-fluid flows}, J. of physics: conference series, 894(2017)012058.

\bibitem{mamprok.izvran18} A.~E.~Mamontov and D.~A.~Prokudin, {\it Solubility of unsteady equations of multi-component viscous compressible fluids}, Izvestiya: Mathematics,
{\bfseries 82}:1 (2018), 140--185.

\bibitem{mamprok.jms} A.~E.~Mamontov and D.~A.~Prokudin, {\it Local solvability of initial-boundary value problem for one-dimensional equations
of polytropic flows of viscous compressible multifluids}, J. of Math. Sciences, {\bfseries 231}:2 (2018), 227--242.

\bibitem{mamprok.local} A.~E.~Mamontov and D.~A.~Prokudin, {\it Unique solvability of initial-boundary value problem for one-dimensional equations of polytropic flows of multicomponent viscous compressible fluids}, Siberian Electr. Math. Reports, {\bfseries 15} (2018), 631--649.

\bibitem{mamprok.muchapz} P.~B.~Mucha, M.~Pokorny and E.~Zatorska, {\it Heat-conducting, compressible mixtures with multicomponent diffusion: Construction of a weak solution}, SIAM J. Math. Anal., {\bfseries 47}:5 (2015), 3747--3797.

\bibitem{mamprok.nigm} R.~I.~Nigmatulin, {\it Dynamics of multiphase media, Vol. 1}, Hemisphere, N.Y., 1990.

\bibitem{mamprok.semiprok17} D.~A.~Prokudin, {\it Unique solvability of initial-boundary value problem for a model system of equations for the polytropic motion of a mixture of viscous compressible fluids}, Siberian Electr. Math. Reports, {\bfseries 14} (2017), 568--585 (in Russian).

\bibitem{mamprok.jfcsprok} D.~A.~Prokudin, {\it Global solvability of the initial boundary value problem for a model system of one-dimensional equations of polytropic flows of viscous compressible fluid mixtures}, J. of Physics: Conference Series, 894 (2017) 012076, 6 P.

\bibitem{mamprok.jaimprok} D.~A.~Prokudin and M.~V.~Krayushkina, {\it Solvability of a stationary boundary value problem for a model system of the equations of barotropic motion of a mixture of compressible viscous fluids}, Journal of Applied and Industrial Mathematics, {\bfseries 10}:3 (2016), 417--428.

\bibitem{mamprok.raj} K.~L.~Rajagopal and L.~Tao, {\it Mechanics of mixtures, Series on Advances in Mathematics for Applied Sciences}, {\bfseries 35},
World Scientific, River Edge, NJ, 1995.

\bibitem{mamprok.weigant92} W.~Weigant, {\it Non-homogeneous boundary value problems for the Navier---Stokes equations of a viscous gas}. Ph. D. Thesis. Barnaul. 1992 (in Russian).
\end{thebibliography}
\end{document}